\documentclass[11pt,oneside]{amsart}
\usepackage[mathscr]{eucal}
\usepackage{amssymb, amsmath, array, amscd}
\usepackage[dvips]{graphicx}
\usepackage{color}
\usepackage{epsfig}
\newtheorem{theorem}{Theorem}
\newtheorem{corollary}[theorem]{Corollary}

\begin{document}

\title[Metrics on Kato surfaces]{Locally conformally K\"ahler metrics\\
on Kato surfaces}

\author{Marco Brunella}

\address{Marco Brunella, Institut de Math\'ematiques de Bourgogne
-- UMR 5584 -- 9 Avenue Savary, 21078 Dijon, France}

\begin{abstract}
We show that every Kato surface admits a locally conformally
K\"ahler metric.
\end{abstract}

\maketitle

In this note we shall prove the following theorem, which improves
and completes our previous result \cite{Bru}.

\begin{theorem} \label{thm}
Every Kato surface admits a locally conformally K\"ahler metric.
\end{theorem}

We refer to \cite{Bel} and \cite{Bru}, and references therein, for some
background on locally conformally K\"ahler metrics on complex surfaces.
We just recall the following fact, which will be used below: a complex
surface $S$ admits a locally conformally K\"ahler metric if and only if
its universal covering $\widetilde S$ admits a K\"ahler metric $\omega$
which transforms conformally under the action of the deck transformations
($g^*(\omega )=c_g\omega$, $c_g\in{\mathbb R}^+$, for every deck
transformation $g$). In particular, we observe the following consequence
of Theorem \ref{thm}.

\begin{corollary} \label{cor}
The universal covering of a Kato surface is K\"ahlerian.
\end{corollary}

It would be interesting to know if this last statement holds for
{\it every} compact complex surface. Note, however, that there are
compact complex surfaces which do not admit locally conformally
K\"ahler metrics \cite{Bel}, hence a generalization of Corollary
\ref{cor} to any compact complex surface should follow a different
path.

The proof of Theorem \ref{thm} is based on a point of view already
exploited in \cite[\S 2]{Bru}, but we shall not need the stability result
proved in \cite [\S 3]{Bru}.

\subsection*{Kato surfaces}

Let us briefly recall the definition and the construction of Kato
surfaces \cite{Kat} \cite{Dlo}. Let $B_r$ denotes the (open) ball of
radius $r$ in ${\mathbb C}^2$, and set $B=B_1$. Take a sequence of
blow-ups
$$\pi : \widehat{B} \longrightarrow B$$
over the origin $0\in B$, subject to the following constraint: at each step,
we blow-up a point belonging to the exceptional divisor created at the
previous step. Take also a holomorphic (up to the boundary) embedding
$$\sigma : \overline{B} \longrightarrow \widehat{B}$$
which sends the origin to a point belonging to the last exceptional divisor
created by $\pi$. Set
$$W=\widehat{B}\setminus \sigma (\overline{B}).$$
We can glue together the two boundary components of $\partial W =
(\partial\widehat{B})\cup (\sigma (\partial B))$ using the real analytic
CR-diffeomorphism
$$\sigma\circ\pi : \partial\widehat{B}\longrightarrow \sigma (\partial B).$$
The result is a minimal compact complex surface $S=S_{\pi ,\sigma}$,
with first Betti number equal to 1 and second Betti number equal to
the number of blow-ups in $\pi$. It is called a {\it Kato surface},
or a {\it surface with a global spherical shell}. Of course, such a
surface cannot be K\"ahlerian, for its first Betti number is odd.

The following remark is very simple, yet of capital importance for
our construction. For every $r\in (0,1]$, the embedding $\sigma$
sends $\overline{B}_r$ into $\widehat{B}_r = \pi^{-1}({B}_r)$:
indeed, the composite map $\pi\circ\sigma : \overline{B}\to {B}$
fixes the origin, and hence by Schwarz Lemma it maps
$\overline{B}_r$ inside $B_r$, for every $r\le 1$, which precisely
means that $\sigma (\overline{B}_r)\subset\pi^{-1}(B_r)$. If we set
$$W_r = \widehat{B}_r \setminus \sigma (\overline{B}_r),$$
then we can glue, as before, the two components of $\partial W_r$
using $\sigma\circ\pi$. The result is the {\it same} compact complex
surface $S$. Indeed, the open subsets $W_r$, $r\in (0,1]$, may be
understood as different fundamental domains in the same universal
covering $\widetilde S$ for the same action of the deck
transformations. This phenomenon is also at the origin of the fact
that the classification of Kato surfaces can be reduced to the
classification of {\it germs} of contracting mappings of the type
$\pi\circ\sigma$ \cite{Dlo}. In fact, we may replace the balls $B_r$
with any domain $U\subset B$ containing $0$ and such that $\sigma
(\overline{U})\subset\pi^{-1}(U)$.

\subsection*{A K\"ahler metric on $\widehat{B}$}

On $\widehat{B}$, there obviously exists a K\"ahler form $\omega_0$, smooth up
to the boundary. In order to prove Theorem \ref{thm}, we firstly modify $\omega_0$
on a neighbourhood of $\sigma (0)$.

Consider $\sigma^*(\omega_0)$. It is a K\"ahler form on $B$, which admits a potential:
$$\sigma^*(\omega_0)={\rm dd^c}\varphi$$
with $\varphi$ smooth and strictly plurisubharmonic on $B$. Up to adding an affine
(hence pluriharmonic) function, we may assume that $\varphi (0)=0$ and
$({\rm d}\varphi )(0)=0$. The Taylor expansion of $\varphi$ at $0$ has the form
$$\varphi (z) = L(z,\overline z) + {\rm Re}M(z,z) + O(\left\| z \right\|^3)$$
where $L$ is the Levi form at $0$ and $M$ is a homogeneous polynomial of degree $2$.
Because ${\rm Re}M$ is pluriharmonic, we may replace $\varphi$ with
$\varphi - {\rm Re}M$, and we obtain in this way a new potential for
$\sigma^*(\omega_0)$, still denoted by $\varphi$, whose Taylor expansion at $0$
is $L(z,\overline z) + O(\left\| z \right\|^3)$. This means that $\varphi$ has
a {\it minimum} point at $0$, of Morse type.

We can choose $\lambda > 0$ and $0<r<r'<1$ such that the function
$$\rho (z) = \lambda (1+\left\| z \right\|^2)$$
is greater than $\varphi$ on $\overline{B}_r$ and smaller than $\varphi$ on
$\partial B_{r'}$. The function $\widetilde\varphi$ on $B$ defined by
$\widetilde\varphi = \varphi$ on $B\setminus B_{r'}$ and $\widetilde\varphi =$
regularized maximum between $\varphi$ and $\rho$ on $B_{r'}$ is therefore smooth,
strictly plurisubharmonic, and equal to $\rho$ on some neighbourhood of
$\overline{B}_r$.

The K\"ahler form $\sigma_*({\rm dd^c}\widetilde\varphi )$ on
$\sigma (B)$ glues to $\omega_0$ outside $\sigma (B)$, giving a new
K\"ahler form $\omega_1$ on $\widehat{B}$. By construction, this
K\"ahler form has the property:
$$\sigma^*(\omega_1) = {\rm dd^c}\rho \quad {\rm on}\ \overline{B}_r.$$

\subsection*{Another K\"ahler metric on $\widehat{B}_r$}

We now modify the previous $\omega_1$ on a neighbourhood of $\partial \widehat{B}_r$.
Choose $s<r$ sufficiently close to $r$, so that $\widehat{B}_s$ still contains
$\sigma (\overline{B}_r)$, i.e. $(\pi\circ\sigma )(\overline{B}_r)\subset B_s$.

\begin{figure}[h]
\centering
\includegraphics[width=11cm,height=6cm]{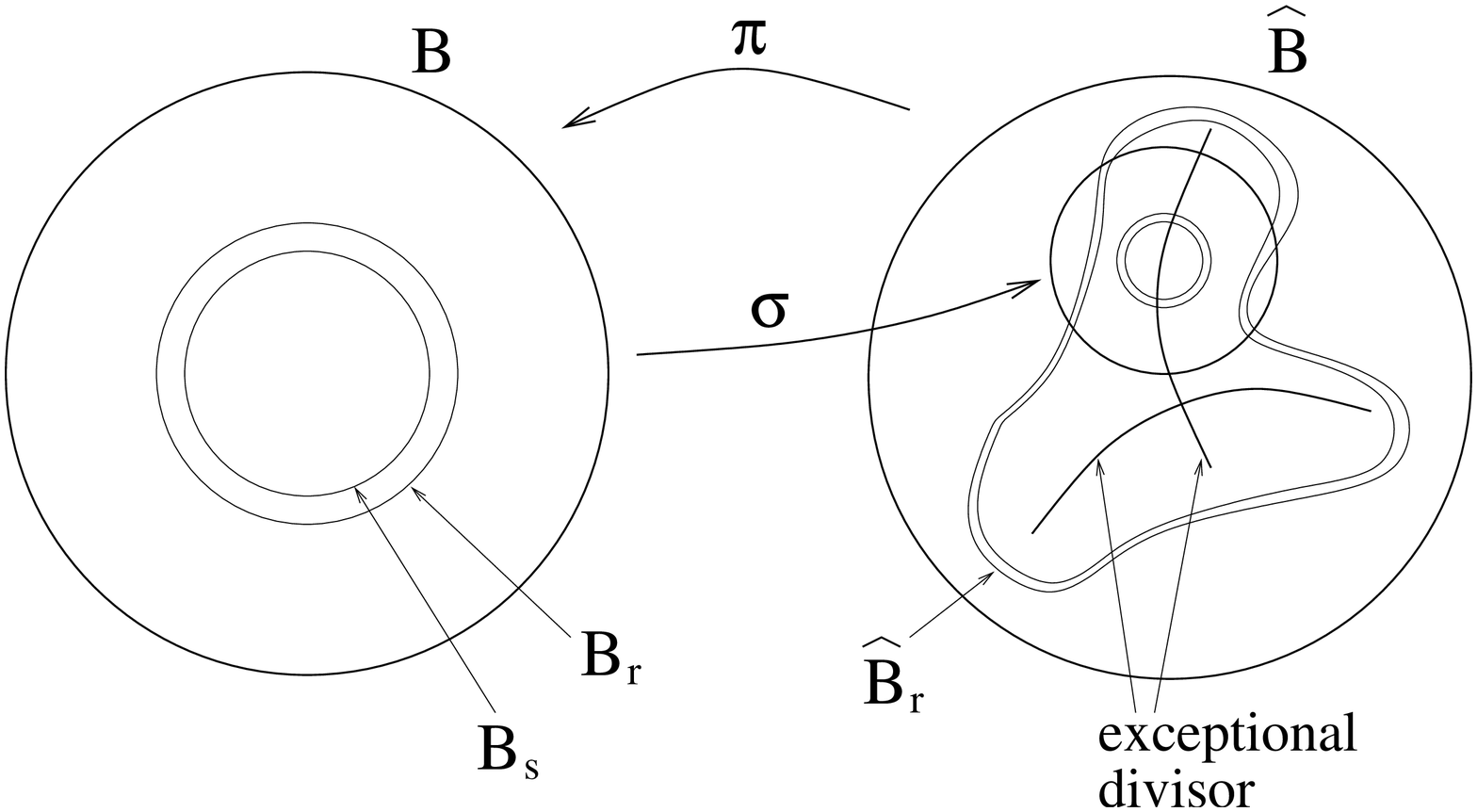}
\end{figure}

Take a potential $\psi$ for the direct image $\pi_*(\omega_1)$:
$$\pi_*(\omega_1) = {\rm dd^c}\psi \quad {\rm on}\ B.$$
Note that $\psi$ is necessarily singular at $0$, but it doesn't matter.
We can choose $c_1\in{\mathbb R}^+$ and $c_2\in{\mathbb R}$ such that the
function $$c_1\rho + c_2$$ is greater than $\psi$ on $\partial B_r$ and smaller
than $\psi$ on $\partial B_s$.
The function $\widetilde\psi$ on $\overline{B}_r$ defined by $\widetilde\psi =\psi$
on $B_s$ and $\widetilde\psi =$ regularized maximum between $\psi$ and $c_1\rho + c_2$
on $\overline{B}_r\setminus B_s$ is therefore smooth, strictly plurisubharmonic,
equal to $c_1\rho + c_2$ on some neighbourhood of $\partial B_r$ and equal to
$\psi$ on some neighbourhood of $\overline{B}_s$.

By pulling-back $\widetilde\psi$ with $\pi$ and taking ${\rm dd^c}$,
we get on $\widehat{B}_r$ a new K\"ahler form $\omega_2$ with the
following two properties:
\begin{enumerate}
\item[(1)] $\omega_2=\omega_1$ on $\widehat{B}_s\supset \sigma (\overline{B}_r)$,
so that $\sigma^* (\omega_2)={\rm dd^c}\rho$ on $\overline{B}_r$;
\item[(2)] $\pi_*(\omega_2)= {\rm dd^c}(c_1\rho + c_2) = c_1{\rm dd^c}\rho$ on
some neighbourhood of $\partial B_r$.
\end{enumerate}

\subsection*{Conformal glueing}

We can now complete the proof of Theorem \ref{thm}. The glueing map
$\sigma\circ\pi : \widehat{B}_r \to \sigma (B_r)$ satisfies the conformal
property
$$(\sigma\circ\pi )^*(\omega_2\vert_U)= \frac{1}{c_1} \omega_2\vert_V$$
where $U\subset \sigma (B_r)$ and $V\subset \widehat{B}_r$ are suitable
neighbourhoods of $\sigma (\partial B_r)$ and $\partial\widehat{B}_r$.
Because the Kato surface $S$ is obtained by identifying these two
neighbourhoods via $\sigma\circ\pi$, we get on $S$ a locally conformally
K\"ahler metric.

\subsection*{Remark}

Because $c_1$ can be chosen arbitrarily large, we see also that, in
the terminology of \cite[Rem. 9]{Bru}, the open subset $I(S)\subset
(0,1)$ contains an interval $(0,\varepsilon )$ (with $\varepsilon$,
in principle, depending on $S$). We do not know if $I(S)$ {\it is}
an interval.

\subsection*{Remark}

The previous construction can be done parametrically: if $\{
S_t\}_{t\in{\mathbb D}}$ is a family of Kato surfaces, then, up to
restricting the base ${\mathbb D}$, we can find a family of locally
conformally K\"ahler metrics $\{\omega_t\}_{t\in{\mathbb D}}$,
smoothly depending on $t$. This is, however, slightly weaker than
the stability result of \cite{Bru}, which states that {\it every}
locally conformally K\"ahler metric $\omega_0$ on $S_0$ extends to a
family $\{\omega_t\}_{t\in{\mathbb D}}$.

\end{document}